\documentclass{elsart3-1}

 \usepackage{graphicx}

\usepackage{amssymb}

\usepackage[english,francais]{babel}
\usepackage[T1]{fontenc}
\usepackage[latin1]{inputenc}
\newtheorem{e-proposition}[theorem]{Proposition}

\newtheorem{e-definition}[theorem]{Definition\rm}

\newtheorem{theoreme}{Th\'eor\`eme}[section]
\newtheorem{lemme}[theoreme]{Lemme}

\renewcommand{\theequation}{\arabic{equation}}
\def\og{\leavevmode\raise.3ex\hbox{$\scriptscriptstyle\langle\!\langle$~}}
\def\fg{\leavevmode\raise.3ex\hbox{~$\!\scriptscriptstyle\,\rangle\!\rangle$}}

\begin{document}
% Vous pouvez mettre dans la prochain ligne la rubrique choisie
% (si vous la connaissez) - meme deux, format : Rubrique1/Rubrique2
\centerline{Statistique}
\begin{frontmatter}

% Titre, auteurs et adresses

% utiliser la commande \thanksref dans \title, \author ou \address
%     pour les notes en bas de page ;

% utiliser la commande \ead pour l'adresse e-mail de chaque auteur
%    (aprËs la commande \auteur) ;

%\title{Title\thanksref{label1}}
%\thanks[label1]{}
% \author{Name\thanksref{label2}}
% \ead{email address}
%
% \thanks[label2]{}
% \address{Address\thanksref{label3}}
% \thanks[label3]{}
\selectlanguage{francais}
\title{Estimation consistante de l'architecture des perceptrons multicouches}

\vspace{-2.6cm}
\selectlanguage{english}
\title{Consistent estimation of the architecture of multilayer perceptrons}

% utiliser les Ètiquettes pour indiquer l'adresse de chaque auteur,
%     s'il y a plusieurs adresses

% \author[label1,label2]{}
% \address[label1]{}
% \address[label2]{}

\author[authorlabel1]{Joseph Rynkiewicz}
\ead{joseph.rynkiewicz@univ-paris1.fr}
\address[authorlabel1]{SAMOS/MATISSE, Universit\'e de Paris-I, 90, rue de Tolbiac 75013 Paris, France, Tél. et Fax : 0144078705}

\begin{abstract}
% Texte de l'abstract en anglais
We consider regression models involving multilayer perceptrons (MLP) with one hidden layer and a Gaussian noise. The estimation of the parameters of the MLP can be done by maximizing the likelihood of the model. In this framework,  it is  difficult to determine the true number of hidden units because the information matrix of Fisher is not invertible if this number is overestimated. However, if the parameters of the MLP are in a compact set, we prove that the minimization of a suitable information criteria leads to consistent estimation of the true number of hidden units. {\it To cite this article:}

\vskip 0.5\baselineskip

\selectlanguage{francais}
% Texte du rÈsumÈ en franÁais
\noindent{\bf R\'esum\'e}
\vskip 0.5\baselineskip
\noindent
On considère des modèles de régression impliquant des perceptrons multicouches (MLP) avec une couche cachée et un bruit gaussien. L'estimation des paramètres du MLP peut être faite en maximisant la vraisemblance du modèle. Dans ce cadre, il est difficile de déterminer le vrai nombre d'unités cachées parce que la matrice d'information de Fisher n'est pas inversible si ce nombre est surestimé. Cependant, si les paramètres du MLP sont dans un ensemble compact, nous prouvons que la minimisation d'un critère d'information convenable permet l'estimation consistante du vrai nombre d'unités cachées.  {\it Pour citer cet article~:}

\end{abstract}
\end{frontmatter}

% Maintenant la version abrÈgÈe en anglais, si prÈsente
%\selectlanguage{english}
%\section*{Abridged English version}
\selectlanguage{francais}
\section{Introduction}
On étudie le comportement asymptotique pour l'estimateur du maximum de vraisemblance d'un modèle de régression utilisant un MLP. On suppose ici qu'il existe un vrai modèle MLP qui a généré les observations.  Lorsque le nombre d'unités cachées du MLP est surestimé, le vrai paramètre du modèle n'est plus identifiable, même à une permutation près. Si les paramètres du MLP ne sont pas bornés à priori,  Fukumizu \cite{Fukumizu} a montré que la statistique du rapport de vraisemblance tendait vers l'infini. Cependant, il est courant de supposer que les paramètres du modèle sont bornés. Dans ce cadre et sous de bonne hypothèses, nous montrons qu'un critère d'information convenablement choisi, par exemple le BIC, est consistant. 

Définissons maintenant notre modèle. Soit les vecteurs de ${\mathbb R}^d$ : $x=(x_1,\cdots,x_d)^T$ et $w_i:=\left(w_{i1},\cdots,w_{id}\right)^T$. La fonction représentée par un MLP avec $k$ unités cachées s'écrit :
\[
F_\theta(x)=\beta+\sum_{i=1}^k a_i\phi\left(b_i+w_i^Tx\right)\\
\]
où $\phi$ est la fonction de transfert qui sera supposée dans toute la suite bornée et trois fois dérivable. On supposera aussi que les dérivées premières, secondes et troisièmes de $\phi$, notées respectivement $\phi^{'}$, $\phi^{''}$ et $\phi^{'''}$, seront bornées. Soit
$\theta=\left(\beta,a_1,\cdots,a_k,b_1,\cdots,b_k,w_{11},\cdots,w_{1d},\cdots,w_{kd}\right)\subset {\mathbb R}^{2k+1+k\times d}$ le vecteur paramètre du modèle. Montrons que si on surestime le nombre d'unités cachées, le vrai paramètre n'est plus identifiable. Supposons, par exemple, que la vraie fonction soit donnée par un MLP avec une seule unité cachée : $F_{\theta^0}(x)=a^0_1\tanh(w^0_{11}x)$ avec $x$ réel et $tanh$ la fonction tangente hyperbolique. Alors, tout paramètre $\theta$ de l'ensemble 
\[
\left\{\theta=\left(w_{11}=w_{21}=w^0_{11}, a_1+a_2=a^0_1, \beta=b_1=b_2=0\right)\right\}
\]
réalisera la fonction $F_{\theta^0}$. Une autre difficulté apparaît lorsque qu'il existe un $w_i$ nul, car la fonction $\phi(b_i+w_i^Tx)$ est alors constante comme $\beta$. Pour éviter ce problème, on restreindra $\Theta$ à l'ensemble des paramètres tels qu'il existe un $\eta$ vérifiant  $\Vert w_i \Vert\geq \eta$, pour tout $w_i\in \Theta$.

On considère une suite de variables aléatoires i.i.d. $Z_i=(X_i,Y_i)$ où $X_i$ a pour loi $q(x)\lambda_d(x)$ avec $\lambda_d$ la mesure de Lebesgue sur ${\mathbb R}^d$ et $q(x)>0$ pour tout $x\in{\mathbb R}^d$. La vraisemblance de l'observation $z:=(x,y)$ s'écrit alors : 
\[
f_\theta(z)=\frac{1}{\sqrt{2\pi\sigma^2}}e^{-\frac{1}{2\sigma^2}\left(y-F_\theta(x)\right)^2}q(x)
\] 
Par souci de simplicité et de concision, on supposera que la variance du bruit $\sigma^2$ est connue. On suppose de plus, que le vrai modèle a, au plus, $M$ unités cachées. L'ensemble des paramètres considérés  est alors noté \( \Theta:=\cup_{1\leq k\leq M} \Theta_k\) avec, pour tout $k$ et un $\eta>0$,
\[
 \Theta_k:=\left\{\theta=\left(\beta,a_1,\cdots,a_k,b_1,\cdots,b_k,w_{11},\cdots,w_{1d},\cdots,w_{kd}\right),\ \forall 1\leq i\leq k, \Vert w_i \Vert\geq \eta\right\}\subset {\mathbb R}^{2k+1+k\times d}
\]
un ensemble supposé compact, c'est-à-dire tel que la norme des vecteurs paramètres de $\Theta_k$ soit bornée. On notera $k^0$ le nombre minimal d'unités cachées tel que $F_{\theta^0}\in \Theta_{k^0}$ représente le vrai modèle et $f(z):=f_{\theta^0}(z)$ la vraie densité des observations. 

\section{Identification de l'architecture du MLP}
Notons $l_n(\theta):=\sum_{i=1}^n\log(f_\theta(z_i))$, on définit l'estimateur du  maximum de vraisemblance pénalisé de $k_0$, comme étant le nombre d'unités cachées $\hat k$ qui maximise 
\(
T_n(k):=\max\{l_n(\theta) : \theta\in\Theta_k\}-p_n(k)
\) 
, où $p_n(k)$ est le terme qui pénalise la log-vraisemblance par le nombre d'unités cachées. On fait maintenant les hypothèses suivantes : 
\begin{description}
\item{H-1 : } les fonctions MLP sont identifiables au sens faible suivant  :
\[
\forall x,\ \beta^0+\sum_{i=1}^{k^0} a^0_i\phi\left(b^0_i+{w^0_i}^Tx\right)=\beta+\sum_{i=1}^k a_i\phi\left(b_i+w_i^Tx\right)\Leftrightarrow \beta=\beta^0\mbox{ et }\sum_{i=1}^{k^0}a^0_i\delta_{(b^0_i,w^0_i)}=\sum_{i=1}^{k}a_i\delta_{(b_i,w_i)}.
\]
où $\delta_x$ est la fonction qui vaut $1$ en $x$ et $0$ partout ailleurs. 
\item{H-2 :} $X$ admet un moment d'ordre 6.
\item{H-3 :} les fonctions de l'ensemble
\[
\begin{array}{l}
\left(\left(x_kx_l\phi^{''}(b^0_i+{w^0_i}^Tx)\right)_{1\leq l \leq k\leq d,\ 1\leq i\leq k^0},\phi^{''}(b^0_i+{w^0_i}^Tx)_{1\leq i\leq k^0}, \right.\\
\left. \left(x_k\phi^{'}(b^0_i+{w^0_i}^Tx)\right)_{1\leq k\leq d,\ 1\leq i\leq k^0}, \left(\phi^{'}(b^0_i+{w^0_i}^Tx)\right)_{1\leq i\leq k^0}\right)
\end{array}
\]
sont linéairement indépendantes dans l'espace de Hilbert $L^2(q\lambda_{d})$. 
\item{H-4 :} $p_n(.)$ est croissante,  $p_n(k_1)-p_n(k_2)\stackrel{n\rightarrow\infty}{\longrightarrow} \infty$ pour tout $k_1>k_2$ et $\lim_{n\rightarrow\infty}\frac{p_n(k)}{n}=0$
\end{description}
On aura alors le résultat suivant : 
\begin{theoreme}
Sous H-1, H-2, H-3 et H-4 :   $\hat k\stackrel{P}{\rightarrow}k_0$. 
\end{theoreme}
\paragraph*{Preuve}
Considérons les fonctions  :
\[
s_\theta(z):=\frac{\frac{f_\theta}{f}(z)-1}{\Vert\frac{f_\theta}{f}-1\Vert_2} 
\mbox{ où } \Vert .\Vert_2 \mbox{ est la norme de } L^2\left(f\lambda_{d+1}\right)
\]
Pour démontrer le théorème, il suffit de montrer que l'ensemble ${\mathbb S}:=\{s_\theta,\ \theta\in \Theta\}$ est une classe de Donsker (cf van der Vaart \cite{Vaart}) et le résultat découlera du théorème 2.1 de Gassiat \cite{Gassiat}. 

Le cas difficile est pour $k\geq k_0$. Nous allons reparamétriser le modèle en utilisant une méthode similaire à celle de Liu et Shao \cite{Shao} pour les modèles de mélange. 
Lorsque $\frac{f_\theta}{f}-1=0$ on a $\beta=\beta^0$ et il existe un vecteur $t=(t_i)_{1\leq i\leq k^0}$ tel que $0=t_0<t_1<\cdots<t_{k^0}\leq k$ et à une permutation près :
\(b_{t_{i-1}+1}=\cdots=b_{t_i}=b^0_i\), \(w_{t_{i-1}+1}=\cdots=w_{t_i}=w^0_i\), \(\sum_{j=t_{i-1}+1}^{t_i}a_j=a_i^0\) et \(a_j=0\) pour $t_{k^0}+1\leq j\leq k$. Définissons $s_i=\sum_{j=t_{i-1}+1}^{t_i}a_j-a_i^0$ et $q_j=\frac{a_j}{\sum_{t_{i-1}+1}^{t_i}a_j}$, on aura alors la reparamétrisation $\theta=\left(\Phi_t,\psi_t\right)$ avec \(\Phi_t=\left(\beta,(b_j)_{j=1}^{t_{k^0}},(w_j)_{j=1}^{t_{k^0}},(s_i)_{i=1}^{k^0},(a_j)_{j=t_{k^0}+1}^{k}\right)\), \(\psi_t=\left((q_j)_{j=1}^{t_{k^0}},(b_j)_{t_{k^0}+1}^{k},(w_j)_{t_{k^0}+1}^{k}\right)\). L'intérêt de cette paramétrisation est que, pour $t$ fixé,  $\Phi_t$ est un paramètre identifiable et toute la non-identifiabilité du modèle sera regroupée dans $\psi_t$. Ainsi $F_{(\Phi_t^0,\psi_t)}$ sera égale à $F_{\theta^0}$ si et seulement si
\[
\begin{array}{ccccccccc}
\Phi^0_t=(\beta^0,&\underbrace{b_1^0,\cdots,b_1^0}&,\cdots,&\underbrace{b_{k^0}^0,\cdots,b_{k^0}^0},&\underbrace{w_1^0,\cdots,w_1^0}&,\cdots,&\underbrace{w_{k^0}^0,\cdots,w_{k^0}^0}&,\underbrace{0,\cdots,0}&,\underbrace{0,\cdots,0})\\
 &t_1& &t_{k^0}-t_{k^0-1}&t_1& &t_{k^0}-t_{k^0-1}&k^0&k-t_{k^0}
\end{array}
\]
On aura alors $\frac{f_\theta}{f}(z)$ qui vaudra
\[
\frac{exp\left(-\frac{1}{2\sigma^2}\left(y-\left(\beta+\sum_{i=1}^{k^0}(s_i+a^0_i)\sum_{j=t_{i-1}+1}^{t_i}q_j\phi(b_j+w_j^Tx)+\sum_{j=t_{k^0}+1}^ka_j\phi(b_j+w_j^Tx)\right)\right)^2\right)}{exp\left(-\frac{1}{2\sigma^2}\left(y-\left(\beta^0+\sum_{i=1}^{k^0} a^0_i\phi(b^0_i+{w^0_i}^Tx)\right)\right)^2\right)}
\] 
\begin{lemme}
Notons $D(\Phi_t,\psi_t):=\Vert\frac{f_{(\Phi_t,\psi_t)}}{f}-1\Vert_2$ et $e(z):=\frac{1}{\sigma^2}\left(y-\left(\beta^0+\sum_{i=1}^{k^0} a^0_i\phi(b^0_i+{w^0_i}^Tx)\right)\right)$ on a alors  l'approximation suivante :
\[
\frac{f_\theta}{f}(z)=1+(\Phi_t-\Phi^0_t)^Tf^{'}_{(\Phi^0_t,\psi_t)}(z)+0.5(\Phi_t-\Phi^0_t)^Tf^{''}_{(\Phi^0_t,\psi_t)}(z)(\Phi_t-\Phi^0_t)+o(D(\Phi_t,\psi_t))
\]
avec 
\[
\begin{array}{l}
(\Phi_t-\Phi^0_t)^Tf^{'}_{(\Phi^0_t,\psi_t)}(z)=
\left(\beta-\beta^0+\sum_{i=1}^{k^0}s_i\phi(b^0_i+{w^0_i}^Tx)+\sum_{i=1}^{k^0}\sum_{j=t_{i-1}+1}^{t_i}q_j\left(b_j-b^0_i\right)a^0_i\phi^{'}(b^0_i+{w^0_i}^Tx)\right.\\
\left.+\sum_{i=1}^{k^0}\sum_{j=t_{i-1}+1}^{t_i}q_j\left(w_{j}-w^0_{i}\right)^Txa^0_i\phi^{'}(b^0_i+{w^0_i}^Tx)+\sum_{j=t_{k^0}+1}^ka_j\phi(b_j+w_j^Tx)\right)e(z)
\end{array}
\]
et
\[
\begin{array}{l}
(\Phi_t-\Phi^0_t)^Tf^{''}_{(\Phi^0_t,\psi_t)}(z)(\Phi_t-\Phi^0_t)=\left(1-\frac{1}{e^2(z)}\right)\left((\Phi_t-\Phi^0_t)^Tf^{'}_{(\Phi^0_t,\psi_t)}(z){f^{'}_{(\Phi^0_t,\psi_t)}}^T(z)(\Phi_t-\Phi^0_t)\right)+e(z)\times\\
\left(\sum_{i=1}^{k^0}\sum_{j=t_{i-1}+1}^{t_i}q_j(b_j-b^0_i)^2a^0_i\phi^{''}(b^0_i+{w^0_i}^Tx)+\sum_{i=1}^{k^0}\sum_{j=t_{i-1}+1}^{t_i}q_j(w_{j}-w^0_{i})^Txx^T(w_{j}-w^0_{i})a^0_i\phi^{''}(b^0_i+{w^0_i}^Tx)\right.\\
\left. +\sum_{i=1}^{k^0}\sum_{j=t_{i-1}+1}^{t_i}(q_jb_j-b^0_i)s_{i}\phi^{'}(b^0_i+{w^0_i}^Tx)+\sum_{i=1}^{k^0}\sum_{j=t_{i-1}+1}^{t_i}(q_jw_{j}-w^0_{i})^Txs_{i}\phi^{'}(b^0_i+{w^0_i}^Tx) \right)
\end{array}
\]
\end{lemme}

Disposant de ce développement asymptotique, exactement de la même façon que dans la preuve de la proposition 3.1 de Dacunha-Castelle et Gassiat \cite{Dacunha} ou bien celle du théorème 4.1 de Kéribin \cite{Keribin}, on montre que le nombre $N(\epsilon)$ d'$\epsilon$-brackets (cf van der Vaart \cite{Vaart}) nécessaire pour recouvrir  $\left\{S_\theta, \theta\in\Theta_k\right\}$ est de l'ordre de $O\left(\frac{1}{\epsilon}^{2k+1+k\times d}\right)$. Cela montre que $\mathbb S$ est une classe de Donsker $\blacksquare$

Sussmann \cite{Sussmann} a montré que si les fonctions $\phi$ sont des fonctions sigmoïdes et si on contraint les paramètres $b_i$ à être positifs   pour tout $1\leq i\leq k$, ceci afin d'éviter une symétrie sur les signes de $(b_i,w_i)$ et $a_i$, alors l'hypothèse H-1 est vérifiée. De plus, en suivant un raisonnement similaire à Fukimizu \cite{Fukumizu2}, on peut montrer que les fonctions sigmoïdes vérifient l'hypothèse H-3. Ce théorème s'applique donc au cas le plus couramment utilisé en pratique.

%\section*{Remerciements}
% Remerciements - texte ici


\begin{thebibliography}{00}

% \bibitem{label1}
% Texte
%\bibitem{label2}
%

\bibitem{Dacunha}
Dacunha-Castelle D. and Gassiat E.,
Testing the order of a model using locally conic parametrization: Population mixtures and stationary ARMA processes,
Ann. Statist. 27 (1999) 1178-1209.
\bibitem{Fukumizu}
Fukumizu, K.,
Likelihood ratio of unidentifiable models and multilayer neural networks,
Ann. Statist. 31 (2003) 833-851.

\bibitem{Fukumizu2}
Fukumizu, K.,
A regularity condition of the information matrix of a multilayer perceptron network,
Neural Networks, 9 (5) (1996) 871-879.

\bibitem{Gassiat}
Gassiat, E., 
Likelihood ratio inequalities with applications to various mixtures,
Ann. Inst. Henri Poincaré  38 (2002) 897-906.

\bibitem{Keribin}
Keribin, C., Consistent estimation of the order of mixture models, Sankhyä A 62 (1) (2000) 49-66.

\bibitem{Shao}
Liu, X. and Shao, Y., Asymptotics for likelihood ratio tests under loss of identifiability,  Ann. Statist. 31 (2003) 807-832.

\bibitem{Sussmann}
Sussmann, H.J., Uniqueness of the weights for minimal feed-forward nets with a given input-output map,
Neural Networks 5 (1992) 589-593.

\bibitem{Vaart}
van der Vaart, A.,
Asymptotic statistics, Cambridge University Press, Cambridge, 1998.

\end{thebibliography}
\end{document}